\documentclass{amsart}

\usepackage{fancyhdr}
\usepackage{lastpage}
\usepackage{stmaryrd,yhmath}

\newcommand{\bdis}{\begin{displaymath}}
\newcommand{\edis}{\end{displaymath}}
\newcommand{\be}{\begin{equation}}
\newcommand{\ee}{\end{equation}}
\newcommand{\mbb}{\mathbb}
\newcommand{\mcal}{\mathcal}

\newcommand{\vp}{\varphi}

\newcommand{\vth}{\vartheta}

\newcommand{\zf}{\zeta\left(\frac{1}{2}+it\right)}

\theoremstyle{definition}

\theoremstyle{remark}
\newtheorem{remark}[]{Remark}

\newtheorem*{mydef1}{{\bf Theorem}}

\newtheorem*{mydef4}{{\bf Corollary}}

\numberwithin{equation}{section}

\begin{document}

\title{Jacob's ladders and new class of integrals containing product of factors $\zeta^2$}

\author{Jan Moser}

\address{Department of Mathematical Analysis and Numerical Mathematics, Comenius University, Mlynska Dolina M105, 842 48 Bratislava, SLOVAKIA}

\email{jan.mozer@fmph.uniba.sk}

\keywords{Riemann zeta-function}

\begin{abstract}
In this paper we obtain new properties of a signal generated by the Riemann zeta-function on the critical line. At the same time we obtain an
asymptotic formula for a new class of transcendental integrals connected with the Riemann zeta-function
\end{abstract}
\maketitle

\section{Introduction}

\subsection{}

In the paper \cite{4} we have obtained the following formula
\be \label{1.1}
\begin{split}
& \frac 1U\int_T^{T+U}\prod_{k=0}^n\left|\zeta\left(\frac 12+i\vp_1^k(t)\right)\right|^2\sim \\
& \sim\prod_{k=0}^n\frac{1}{\vp_1^k(T+U)-\vp_1^k(T)}\int_{\vp_1^k(T)}^{\vp_1^k(T+U)}\left|\zf\right|^2{\rm d}t, \\
& U\in \left(\left. 0,\frac{T}{\ln^2T}\right]\right.,\ T\to\infty  .
\end{split}
\ee
A motivation for this formula was the well-known multiplicative formula
\bdis
M\left(\prod_{k=1}^n X_k\right)=\prod_{k=1}^n M(X_k)
\edis
from the theory of probability where $X_k$ are the independent random variables and $M$ is the population mean. Some new art of the asymptotic
independence of the partial functions
\bdis
\left|\zf\right|^2,\ t\in \left[\vp_1^k(T),\vp_1^k(T+U)\right],\ k=0,1,\dots ,n
\edis
is expressed by this formula.

\subsection{}

For example, by using the mean-value theorem in (\ref{1.1}) we obtain
\bdis
\begin{split}
& \left|\zeta\left(\frac 12+i\vp_1^n(\bar{t}_n)\right)\right|^2\frac 1U\int_T^{T+U}\prod_{k=0}^{n-1}
\left|\zeta\left(\frac 12+i\vp_1^k(t)\right)\right|^2{\rm d}t\sim \\
& \sim \frac{1}{\vp_1^n(T+U)-\vp_1^n(T)}\int_{\vp_1^n(T)}^{\vp_1^n(T+U)}\left|\zf\right|^2{\rm d}t \times \\
& \times\prod_{k=0}^{n-1}\frac{1}{\vp_1^k(T+U)-\vp_1^k(T)}\int_{\vp_1^k(T)}^{\vp_1^k(T+U)}\left|\zf\right|^2{\rm d}t
\end{split}
\edis
i. e. (see (\ref{1.1}), $n\mapsto n-1$)
\be \label{1.2}
\left|\zeta\left(\frac 12+i\vp_1^n(\bar{t}_n)\right)\right|^2\sim
\frac{1}{\vp_1^n(T+U)-\vp_1^n(T)}\int_{\vp_1^n(T)}^{\vp_1^n(T+U)}\left|\zf\right|^2{\rm d}t .
\ee
But
\bdis
\left|\zeta\left(\frac 12+i\vp_1^n(\bar{t}_n)\right)\right|^2,\ \bar{t}_n\in (T,T+U)
\edis
is the mean value with respect to the set of functions
\be \label{1.3}
\left\{\left|\zeta\left(\frac 12+i\vp_1^0(t)\right)\right|^2,\dots ,\left|\zeta\left(\frac 12+i\vp_1^{n-1}(t)\right)\right|^2\right\}
\ee
i. e. $\bar{t}_n$ is the nonlinear functional
\be \label{1.4}
\bar{t}_n=f_n[\vp_1^1,\dots ,\vp_1^{n-1}];\ \vp_1^0(t)=t
\ee
that is defined on the continuum set of points
\be \label{1.5}
(\vp_1^1,\dots ,\vp_1^{n-1}) .
\ee
Let us remind that there is the continuum set of the Jacob's ladders (see \cite{1} generating the set of the iterations (\ref{1.5})). At the
same time it follows from (\ref{1.2}) that
\be \label{1.6}
\left|\zeta\left(\frac 12+i\vp_1^n(t_n)\right)\right|^2\sim \left|\zeta\left(\frac 12+i\vp_1^n(\tau)\right)\right|^2,\
\tau\in(T,T+U)
\ee
where $\tau$ is completely independent on the set of points (\ref{1.5}).

\begin{remark}
Thus, the mean-value (\ref{1.2}) with respect to the set of functions (\ref{1.3}) is asymptotically independent on this set.
\end{remark}

\begin{remark}
Now, let $k:\ 1<k<n$. Then we have the mean-value
\be \label{1.7}
\left|\zeta\left(\frac 12+i\vp_1^k(\bar{t}_k)\right)\right|^2
\ee
of the inner factor of the product in (\ref{1.1}) with respect to the set (comp. (\ref{1.3}))
\be \label{1.8}
\begin{split}
& \left\{ \left|\zeta\left(\frac 12+i\vp_1^0(t)\right)\right|^2,\dots ,\left|\zeta\left(\frac 12+i\vp_1^{k-1}(t)\right)\right|^2,
\left|\zeta\left(\frac 12+i\vp_1^{k+1}(t)\right)\right|^2\dots  \right. \\
&
\left. \dots , \left|\zeta\left(\frac 12+i\vp_1^n(t)\right)\right|^2\right\},
\end{split}
\ee
where (comp. (\ref{1.4}))
\bdis
\bar{t}_k=g_k[\vp_1^1,\dots,\vp_1^{k-1},\vp_1^{k+1},\dots,\vp_1^n]
\edis
is the functional defined on the continuum set of points
\bdis
\left(\vp_1^1,\dots,\vp_1^{k-1},\vp_1^{k+1},\dots ,\vp_1^n\right) .
\edis
In this case the mean-value (\ref{1.7}) is not, probable, asymptotically independent on the set (\ref{1.8}).
\end{remark}

\subsection{}

In this paper we obtain new properties of the signal
\bdis
Z(t)=e^{i\vth(t)}\zf,\ \vth(t)=-\frac t2\ln\pi+\text{Im}\ln\Gamma\left(\frac 14+i\frac t2\right)
\edis
generated by the Riemann zeta-function on the critical line. Namely, we obtain an asymptotic formula for a new class of transcendental
integrals of the type
\bdis
\int_T^{T+U}F[\vp_1^{n+1}(t)]\prod_{k=0}^n\left|\zeta\left(\frac 12+i\vp_1^k(t)\right)\right|^2{\rm d}t,\ U\in \left(\left. 0,\frac{T}{\ln^2T}\right]\right. ,
\edis
where
$$F(w),\ w\in [\vp_1^{n+1}(T),\vp_1^{n+1}(T+U)],\ F(w)\geq 0\ (\leq 0)$$
is an arbitrary Lebesgue integrable function. In this direction, we obtain, for example, new asymptotic formulae generalizing our formulae containing
the factors
\bdis
\left|\zeta\left(\frac 12+i\vp_1(t)\right)\right|^4,\ \left\{\arg\zeta\left(\frac 12+i\vp_1(t)\right)\right\}^{2l},\ l\in\mbb{N} .
\edis
Further we obtain the new effect for the macroscopic domains, i. e. for
\bdis
U\in \left[ T^{\frac 13+\epsilon}, \frac{T}{\ln^2T}\right] .
\edis
Namely:
\begin{itemize}
\item[(a)] the transformations
\bdis
[T,T+U]\to [\vp_1^k(T),\vp_1^k(T+U)],\ k=1,\dots ,n+1
\edis
asymptotically preserve the measure (the length) of the segment $[T,T+U]$, i. e. that
\bdis
|[\vp_1^k(T),\vp_1^k(T+U)]|\sim U,\ k=1,\dots ,n+1,\ T\to\infty ;
\edis
\item[(b)] the segments
\bdis
\{[\vp_1^k(T),\vp_1^k(T+U)]\}_{k=0}^{n+1}
\edis
are distributed with an exact asymptotic regularity.
\end{itemize}

\section{The result}

\subsection{}

Let us remind that the formula (comp. \cite{4}, (3.7), (3.8))
\be \label{2.1}
\tilde{Z}^2(t)=\frac{{\rm d}\vp_1(t)}{{\rm d}t},\ \vp_1(t)=\frac 12\vp(t),\ t\geq T_0[\vp]
\ee
and
\be \label{2.2}
\tilde{Z}^2(t)=\frac{Z^2(t)}{2\Phi'_\vp[\vp(t)]}=\frac{\left|\zf\right|^2}{\left\{1+\mcal{O}\left(\frac{\ln\ln t}{\ln t}\right)\right\}\ln t}
\ee
where $\vp(t)$ is the Jacob's ladder, i. e. the exact solution of the nonlinear integral equation
\be \label{2.3}
\begin{split}
& \int_0^{\mu[x(T)]} Z^2(t)e^{-\frac{2}{x(T)}}{\rm d}t=\int_0^T Z^2(t){\rm d}t , \\
& \mu(y)\geq 7y\ln y,\ \mu(y)\to y=\vp_\mu(T)=\vp(T)
\end{split}
\ee
(see \cite{1}). Next, we have (see \cite{4}, (2.1))
\bdis
\begin{split}
& y=\frac 12\vp(t)=\vp_1(t);\ \vp_1^0(t)=t,\ \vp_1^1(t)=\vp_1(t), \\
& \vp_1^2(t)=\vp_1(\vp_1(t)),\dots ,\vp_1^k(t)=\vp_1(\vp_1^{k-1}(t)),\dots
\end{split}
\edis
where $\vp_1^k(t)$ stands for the $k$th iteration of the Jacob's ladder
\bdis
y=\vp_1(t)
\edis
(of course, $\vp_1^k(t),\ t\in [T,T+U]$ are the increasing functions.) The following Theorem holds true.

\begin{mydef1}
Let
\be \label{2.4}
U\in \left(\left. 0,\frac{T}{\ln^2T}\right]\right. .
\ee
Then for every fixed $n\in\mbb{N}$ and for every Lebesgue-integrable function
\bdis
F(t),\ t\in [\vp_1^{n+1}(T),\vp_1^{n+1}(T+U)],\quad F(t)\geq 0 \ (\leq 0)
\edis
we have
\be \label{2.5}
\begin{split}
& \int_T^{T+U}F[\vp_1^{n+1}(t)]\prod_{k=0}^n\left|\zeta\left(\frac 12+i\vp_1^k(t)\right)\right|^2{\rm d}t\sim \\
& \sim \left\{\int_{\vp_1^{n+1}(T)}^{\vp_1^{n+1}(T+U)}F(t){\rm d}t\right\}\ln^{n+1}T,\ T\to \infty
\end{split}
\ee
where
\be \label{2.6}
\vp_1^k(T+U)-\vp_1^k(T)<\frac{1}{2n+3}\frac{T}{\ln T},\ k=1,\dots ,n+1 ,
\ee
\be \label{2.7}
\vp_1^k(T+U)-\vp_1^k(T)>0.2\times\frac{T}{\ln T},\ k=0,1,\dots ,n .
\ee
Next, in the macroscopic case (comp. (\ref{2.4}))
\be \label{2.8}
U\in \left[T^{\frac 13+\epsilon},\frac{T}{\ln^2T}\right] ,
\ee
we have more exact information
\be \label{2.9}
|[\vp_1^k(T),\vp_1^k(T+U)]|=\vp_1^k(T+U)-\vp_1^k(T)\sim U,\ k=1,\dots , n ,
\ee
\be \label{2.10}
\vp_1^k(T)-\vp_1^{k+1}(T+U)\sim (1-c)\frac{T}{\ln T},\ k=0,1,\dots ,n ,
\ee
\be \label{2.11}
\begin{split}
& \rho\{[\vp_1^{k-1}(T),\vp_1^{k-1}(T+U)];[\vp_1^{k}(T),\vp_1^{k}(T+U)]\}\sim(1-c)\frac{T}{\ln T}, \\
& k=1,\dots ,n+1
\end{split}
\ee
where $\rho$ denotes the distance of the corresponding segments.
\end{mydef1}

\begin{remark}
In the macroscopic case (\ref{2.8}) the following is true. The system of iterated segments
\be \label{2.12}
[\vp_1^{n+1}(T),\vp_1^{n+1}(T+U)], [\vp_1^{n}(T),\vp_1^{n}(T+U)],\dots ,[T,T+U]
\ee
is the disconnected set and its components are:
\begin{itemize}
\item[(a)] asymptotically equal (see (\ref{2.9})),
\item[(b)] distributed with the asymptotic regularity from the right to the left (see (\ref{2.10}) -- (\ref{2.12})).
\end{itemize}
\end{remark}

\begin{remark}
Every Jacob's ladder
\bdis
\vp_1(t)=\frac 12\vp(t)
\edis
(see (\ref{2.1})) where $\vp(t)$ is the exact solution of the nonlinear integral equation (\ref{2.3}) is the asymptotic solution of the following
nonlinear integro-iteration equation
\be \label{2.13}
\begin{split}
& \frac 1U\int_T^{T+U}F[x^{n+1}(t)]\prod_{k=0}^n\left|\zeta\left(\frac 12+ix^k(t)\right)\right|^2{\rm d}t= \\
& =\left\{\int_{x^{n+1}(T)}^{x^{n+1}(T+U)}F(t){\rm d}t\right\}\ln^{n+1}T
\end{split}
\ee
(comp. (\ref{2.5})) where
\bdis
x_0(t)=t,\ x^1(t)=x(t),\ x^2(t)=x(x(t)),\dots
\edis
i. e. the function $x^k(t)$ is the $k$th iteration of the function $x(t)$, (comp. (\ref{2.13}) with \cite{3}, (11.1), (11.4), (11.6), (11.8),
\cite{4}, (2.5) and \cite{5}, (2.6)).
\end{remark}

\begin{remark}
Similar remarks like Remark 1 -- Remark 2 hold true also when speaking \emph{on the independence of the mean-value}.
\end{remark}

\subsection{}

By (\ref{2.9}) and the formula (see \cite{4}, (3.5))
\bdis
t-\vp_1^{n+1}(t)\sim(1-c)(n+1)\frac{t}{\ln t}
\edis
we obtain easily (for example) from (\ref{2.5}) the following

\begin{mydef4}
In the macroscopic case (\ref{2.8}) we have
\be \label{2.14}
\int_T^{T+U}\prod_{k=0}^n \left|\zeta\left(\frac 12+i\vp_1^k(t)\right)\right|^2{\rm d}t\sim U\ln^{n+1}T,\ T\to\infty ,
\ee
\be \label{2.15}
\begin{split}
& \int_T^{T+U} \left|\zeta\left(\frac 12+i\vp_1^{n+1}(t)\right)\right|^4\prod_{k=0}^n \left|\zeta\left(\frac 12+i\vp_1^k(t)\right)\right|^2{\rm d}t \sim \\
& \sim \frac{1}{2\pi^2}U_1\ln^{n+5}T,\ U_1=T^{\frac 78+\epsilon} ,
\end{split}
\ee
\be \label{2.16}
\begin{split}
& \int_T^{T+U}\left\{\arg\zeta\left(\frac 12+i\vp_1^{n+1}(t)\right)\right\}^{2l}\prod_{k=0}^n \left|\zeta\left(\frac 12+i\vp_1^k(t)\right)\right|^2{\rm d}t\sim \\
& \sim \frac{(2l)!}{l!4^l}U\ln^{n+1}T(\ln\ln T)^l,\ U\in\left[ T^{\frac 12+\epsilon},\frac{T}{\ln^2T}\right] ,
\end{split}
\ee
\be \label{2.17}
\begin{split}
& \int_T^{T+U}\left\{ S_1\left[\vp_1^{n+_1}(t)\right]\right\}^{2l}\prod_{k=0}^n \left|\zeta\left(\frac 12+i\vp_1^k(t)\right)\right|^2{\rm d}t\sim \\
& \sim d_lU\ln^{n+1}T,\ U\in\left[ T^{\frac 12+\epsilon},\frac{T}{\ln^2T}\right] ,
\end{split}
\ee
for every fixed $l,n\in\mbb{N}$ where
\bdis
S_1(T)=\int_0^{T}S(t){\rm d}t,\ S(t)=\frac 1\pi\arg\zf ,
\edis
and the argument is defined by the usual way (comp. \cite{6}, p. 179).
\end{mydef4}

\begin{remark}
The formulae (\ref{2.15}) -- (\ref{2.17}) are can be understand as generalization of our formulae \cite{3}, (8.3), \cite{5}, Lemma 2, (5.4), (5.5). The formula (\ref{2.14}) can be compared with the formula (2.3)
from the paper of reference \cite{4} in the macroscopic case. The small improvements of the Heath-Brown exponent $\frac 78$ in (\ref{2.15}) are irrelevant for our purpose.
\end{remark}

\section{Proof of Theorem}

\subsection{}

By the formula (see \cite{4}, (3.9))
\bdis
\prod_{k=0}^n\tilde{Z}^2[\vp_1^k(t)]=\frac{{\rm d}\vp_1^{n+1}}{{\rm d}t}
\edis
we obtain
\bdis
\begin{split}
& \int_T^{T+U}F[\vp_1^{n+1}(t)]\prod_{k=0}^n\tilde{Z}^2[\vp_1^k(t)]{\rm d}t= \\
& =\int_T^{T+U}F[\vp_1^{n+1}(t)]{\rm d}\vp_1^{n+1}(t)=\int_{\vp_1^{n+1}(T)}^{\vp_1^{n+1}(T+U)}F(t){\rm d}t ,
\end{split}
\edis
i. e.
\be \label{3.1}
\int_T^{T+U}F[\vp_1^{n+1}(t)]\prod_{k=0}^n\tilde{Z}^2[\vp_1^k(t)]{\rm d}t=\int_{\vp_1^{n+1}(T)}^{\vp_1^{n+1}(T+U)}F(t){\rm d}t .
\ee
Since (see \cite{4}, (3.3), (3.6))
\bdis
\begin{split}
& t>\vp_1^1(t)>\vp_1^2(t)>\dots >\vp_1^{n+1}(t) , \\
& (1-\epsilon)T<\vp_1^{n+1}(T)<T
\end{split}
\edis
then
\bdis
(1-\epsilon)T<\vp_1^{n+1}(T)<T+U,\ U\in \left(\left. 0,\frac{T}{\ln^2T}\right]\right. .
\edis
Consequently,
\be \label{3.2}
T'\in (\vp_1^{n+1}(T),T+U) \ \Rightarrow \ \ln T'=\ln T+\mcal{O}(1) .
\ee
Now, if we use the mean-value theorem on the left-hand side of (\ref{3.1}) we obtain (see (\ref{2.2}), (\ref{3.2}))
\be \label{3.3}
\begin{split}
& \int_{T}^{T+U}F[\vp_1^{n+1}(t)]\prod_{k=0}^n\tilde{Z}^2[\vp_1^k(t)]{\rm d}t \sim \\
& \sim \frac{1}{\ln^{n+1}T}\int_T^{T+U}\prod_{k=0}^n\left|\zeta\left(\frac 12+i\vp_1^k(t)\right)\right|^2{\rm d}t  .
\end{split}
\ee
Hence, from (\ref{3.1}) by (\ref{3.3}) the asymptotic formula (\ref{2.5}) follows.

\subsection{}

Let us remind the estimates (see \cite{4}, (3.15))
\bdis
\vp_1^k(T+U)-\vp_1^k(T)<\frac{2k+1}{2n+1}\frac{T}{\ln T}\leq \frac{T}{\ln T},\ k=1,\dots , n .
\edis
It is clear that by the substitution
\bdis
2n+1\to (2n+3)^2
\edis
(for example) in \cite{4}, part 3.4 we obtain the estimates
\be \label{3.4}
\vp_1^k(T+U)-\vp_1^k(T)<\frac{2k+1}{(2n+3)^2}\frac{T}{\ln T}<\frac{1}{2n+3}\frac{T}{\ln T},\ k=1,\dots ,n+1 ,
\ee
i. e. the inequality (\ref{2.6}) holds true. \\

Next we have (see \cite{4}, (3.4))
\be \label{3.5}
\vp_1^k(T)-\vp_1^{k+1}(T+U)+\vp_1^{k+1}(T+U)-\vp_1^{k+1}(T)>(1-\epsilon)(1-c)\frac{T}{\ln T} .
\ee
Consequently we have (see (\ref{3.4}))
\bdis
\begin{split}
& \vp_1^k(T)-\vp_1^{k+1}(T+U)>(1-\epsilon)(1-c)\frac{T}{\ln T}-\{ \vp_1^{k+1}(T+U)-\vp_1^{k+1}(T) \}> \\
& >(1-\epsilon)(1-c)\frac{T}{\ln T}-\frac{1}{2n+3}\frac{T}{\ln T}>\left( 1-c-\frac{1}{2n+3}-\epsilon\right)\frac{T}{\ln T}\geq \\
& \geq \left( 1-c-\frac 15-\epsilon\right)\frac{T}{\ln T}>(0.22-\epsilon)\frac{T}{\ln T}>0.2\frac{T}{\ln T} ,
\end{split}
\edis
since
\bdis
c<0.58 \ \Rightarrow \ 1-c>0.42>\frac 15=0.2 ,
\edis
i. e. the inequality (\ref{2.7}) holds true.

\subsection{}

We use the Hardy-Littlewood-Ingham formula
\be \label{3.6}
\int_T^{T+U} Z^2(t){\rm d}t\sim U\ln T,\ U\in \left[ T^{\frac 13+\epsilon},\frac{T}{\ln^2T}\right]
\ee
where $\frac 13$ is the Balasubramanian exponent (the small improvements of the exponent $\frac 13$ are irrelevant for our purpose), and our formula (see \cite{3}, (2.5))
\be \label{3.7}
\int_T^{T+U}Z^2(t){\rm d}t\sim \{\vp_1(T+U)-\vp_1(T)\}\ln T .
\ee
Comparing the formulae (\ref{3.6}) and (\ref{3.7}) we obtain
\bdis
\vp_1^1(T+U)-\vp_1^1(T)=\vp_1(T+U)-\vp_1(T)\sim U .
\edis
Similarly, by comparison in the cases (see (\ref{2.6}))
\bdis
T\to\vp_1^1(T),\ T+U\to \vp_1^1(T+U); \dots
\edis
we obtain
\be \label{3.8}
\vp_1^k(T+U)-\vp_1^k(T)\sim U,\ k=1,\dots ,n+1
\ee
i. e. the formula (\ref{2.9}) holds true.

\subsection{}

Next, we have (see (\ref{2.4}), (\ref{3.5}), (\ref{3.8}))
\bdis
\begin{split}
& \vp_1^k(T)-\vp_1^{k+1}(T+U)+\vp_1^{k+1}(T+U)-\vp_1^{k+1}(T)\sim(1-c)\frac{T}{\ln T} , \\
& \vp_1^k(T)-\vp_1^{k+1}(T+U)\sim (1-c)\frac{T}{\ln T}-\{\vp_1^{k+1}(T+U)-\vp_1^{k+1}(T)\}\sim \\
& \sim (1-c)\frac{T}{\ln T}-U\sim (1-c)\frac{T}{\ln T}
\end{split}
\edis
i. e. the formula (\ref{2.10}) holds true. The proposition (\ref{2.11}) follows from (\ref{2.10}) .

\thanks{I would like to thank Michal Demetrian for helping me with the electronic version of this work.}

\end{document}